\documentclass[12pt]{article}
\usepackage{amsmath}
\usepackage{amssymb}
\usepackage{dsfont}
\usepackage{cite}
\newsymbol \blackbox 1004
\newcommand{\eh}{\hfill}\newlength{\sperr}

\def\di{\mathrm{diag}}
\def\prt{\partial}
\def\vr{\varrho}

\def\nn{\nonumber}
\def\a{\alpha}
\def\b{\beta}
\def\g{\gamma}

\def\Lam{\Lambda}

\def\om{\omega}
\def\Om{\Omega}

\def\vp{\varphi}
\def\vt{\vartheta}
\def\ve{\varepsilon}
\def\wh{\widehat}
\def\wt{\widetilde}
\def\ov{\overline}

\def\p{\partial}

\def\BC{{\mathbb C}}
\def\BR{{\mathbb R}}
\def\BN{{\mathbb N}}
\def\clp{{\mathcal P}}

\def\clh{{\mathcal H}}

\def\cln{{\mathcal N}}
\def\clq{{\mathcal Q}}

\def\clv{{\mathcal V}}
\def\cld{{\mathcal D}}

\def\M{{\mathfrak M}}

\newcommand{\E}{\mathrm{e}}
\newcommand{\I}{\mathrm{i}}
\def\mf{\mathfrak}

\newtheorem{Pa}{Paper}[section]
\newtheorem{Tm}[Pa]{{\bf Theorem}}

\newtheorem{Cy}[Pa]{{\bf Corollary}}
\newtheorem{Rk}[Pa]{{\bf Remark}}

\newtheorem{Dn}[Pa]{{\bf Definition}}

\newtheorem{Pn}[Pa]{{\bf Proposition}}

\title{Evolution of Weyl
functions and initial-boundary value problems}

\author{A.L. Sakhnovich}

\date{}
\begin{document}
\maketitle

\begin{abstract} This review is dedicated to some recent results on Weyl theory, inverse
problems, evolution of the Weyl functions and  applications to integrable wave equations in a semistrip
and quarter-plane. For overdetermined initial-boundary value problems, we consider some approaches,
which help to reduce the number of  the initial-boundary conditions. The interconnections between dynamical
and spectral Dirac systems, between response and Weyl functions are studied as well.

\end{abstract}


\tableofcontents


\section{Introduction} \label{intro}
\setcounter{equation}{0} 
 Initial-boundary value problems for linear and nonlinear wave equations are interesting, difficult
and have many applications since initial-boundary conditions should be taken into account in the study
of various applied problems connected with wave processes. Weyl theory could be useful in this study (and, vice versa, 
initial-boundary value problems are of interest in Weyl theory).  A number of related results and references
was presented in our book  \cite[Ch. 6]{SaSaR}. Here we mainly review the results that appeared after the
publication of \cite{SaSaR} although some results from \cite{SaSaR} and some earlier results that are not
contained in \cite{SaSaR} are also discussed for the sake of completeness.

 In order to make the review reader-friendly and self-sufficient, we provide short schemes of the proofs
 for some of the main statements and start the review with the section
 "Preliminaries: Weyl functions, inverse problems and zero curvature representation". However,
  even this section contains some quite recent important results.
 
 Zero curvature representation \cite{AKNS, FT, Nov, ZM0}  is  in a certain sense an analogue of the famous Lax pairs. It is a compatibility
 condition for linear systems, which are so called auxiliary systems for nonlinear integrable
 wave equations. This compatibility condition is discussed in Subsection \ref{Comco} of the section "Preliminaries ...".
 Here already, an initial-boundary value problem, which is basic for the study
 of the evolution of Weyl (Weyl-Titchmarsh) functions, appears. 
 
 Section \ref{GW} is dedicated to generalized Weyl functions and corresponding inverse problems.
 In Section \ref{SG}  and in Subsection \ref{unbound} we consider  some applications of   evolution formulas for Weyl  (and generalized Weyl) functions to the problems of uniqueness and
existence and to the unbounded solutions of  wave equations. 

We note that Cauchy problems (in particular, Cauchy problems for integrable nonlinear
wave equations) are studied much more thoroughly  than initial-boundary value problems.
In fact, too many  initial-boundary conditions are usually required for the study of  initial-boundary value problems
(see, e.g., \cite{Ash, BoFo, Fo, SaSaR} and references therein), and the problem becomes overdetermined. 
Thus, a crucial step here is the reduction of  the initial-boundary conditions. Some important approaches to this
reduction are discussed in Section \ref{Reduc}.

Finally, Section \ref{Dyn} is dedicated to the initial-boundary value problem for dynamical Dirac system
and its connections with Weyl theory.

Appendix \ref{AppB} contains some definitions and references from the theory of quasi-analytic functions. 

{\bf Notations.} As usual,  $\BR$ stands for the real axis,
$\BR_+=(0, \, \infty)$,
$\BC$ stands for the complex plain, and
$\BC_+$ for the open upper semi-plane. 
We use notations $\BC_M = \{ z \, \colon \, \Im(z) >  M \}$ and $\BC^{-}_M = \{ z \, \colon \, \Im(z) < - M \}$.
The equality $D=\di\{d_1, ...\}$ means that $D$ is a diagonal matrix with the entries
$d_1, \ldots$ on the main diagonal and $I_m$ denotes the $m\times m$ identity matrix.
"Locally" (e.g., locally summable) with respect to the semiaxis $[0,\infty)$ means: on all the finite intervals $[0,l]$ (e.g.,  summable
on all $[0,l]$). The class of linear bounded operators acting from the normed space $\clh_1$ into the normed space $\clh_2$
is denoted by $B(\clh_1,\clh_2)$ and we write simply $B(\clh)$ when $\clh_1=\clh_2=:\clh$.
We say that a matrix function is boundedly differentiable when its derivative is bounded in the matrix norm,
and an operator is boundedly invertible when the inverse operator exists and is bounded.
The notation ${\mathrm{l.i.m.}}$ stands for the entrywise limit of a matrix function in the norm of  $L^2(0,\ell)$, 
$0<\ell \leq~\infty$.
The notation $\mf M$
stands for
the operator mapping a Weyl function (or generalized Weyl function, depending on the context) $\vp(z)$
on the corresponding potential
(e.g., $\M (\varphi) = v$ for Dirac systems \eqref{p1} and \eqref{s1} where $V$ has the form \eqref{p2}).
That is, $\M (\varphi) $ is the solution of the inverse problem to recover the potential from the   Weyl function.

\section{Preliminaries: Weyl functions, inverse problems and zero curvature representation} \label{Prel}
\setcounter{equation}{0} 
\subsection{Main scheme}\label{Mais}
One of the most important auxiliary linear systems in soliton theory is selfadjoint Dirac (also called ZS or AKNS) system 
\begin{align} &       \label{p1}
\frac{d}{dx}y(x, z )=\I (z j+jV(x))y(x,
z ), \qquad
x \geq 0,
\end{align} 
where
\begin{align} &   \label{p2}
j = \left[
\begin{array}{cc}
I_{m_1} & 0 \\ 0 & -I_{m_2}
\end{array}
\right], \hspace{1em} V= \left[\begin{array}{cc}
0&v\\v^{*}&0\end{array}\right],  \quad m_1+m_2=:m,
 \end{align} 
$I_{m_i}$ is the $m_i \times
m_i$ identity
matrix and $v(x)$ is an $m_1 \times m_2$ matrix function.

For instance, the well-known matrix defocusing nonlinear Schr\"odinger (defocusing NLS or dNLS) equation
\begin{align}  \label{p3} &
2v_t=\I(v_{xx}-2vv^*v) \qquad \Big(v_t:=\frac{\p}{\p t}v\Big)
\end{align}
admits \cite{ZS1, ZS2} (i.e., is equivalent to) zero curvature representation 
\begin{align}  \label{p4} &
G_t-F_x+[G,F]=0 \qquad ([G,F]:=GF-FG),
\\  \label{p5}&
G=\I (zj+jV), \quad
F=-\I\big(z^2 j+ z jV-\big(\I V_x-jV^2\big)/2\big),
\end{align}
where $V$ has the form \eqref{p2} and $v=v(x,t)$.
For this reason, systems $y_x=Gy$ and $y_t=Fy$ are called auxiliary systems for dNLS,
and it  is easy to see that $y_x=Gy$ has the form \eqref{p1} for each fixed $t$.
We note that $V$ is called the potential of system \eqref{p1}. (Sometimes, for convenience, $v$ is also
called the potential.)

{\it Weyl solutions}  are squarely summable on $(0, \infty)$ solutions of  Dirac system \eqref{p1}.
Clearly, these solutions may be presented as linear combinations of  the columns of the 
normalized by the condition
\begin{align} &   \label{p6}
u(0,z)=I_m
\end{align}
fundamental solution $u$ of \eqref{p1}. Thus, Weyl-Titchmarsh $($or simply Weyl$)$ function
of \eqref{p1} is introduced in terms of $u$.
\begin{Dn} \label{DnW} Let Dirac system \eqref{p1} on $[0, \, \infty)$
 be given and assume that $v$ is  locally summable.
Then Weyl function is an $m_2 \times m_1$ holomorphic matrix function, which satisfies the inequality
\begin{align}&      \label{p7}
\int_0^{\infty}
\begin{bmatrix}
I_{m_1} & \vp(z)^*
\end{bmatrix}
u(x,z)^*u(x,z)
\begin{bmatrix}
I_{m_1} \\ \vp(z)
\end{bmatrix}dx< \infty, \quad z \in \BC_+ .
\end{align}
\end{Dn}
The following proposition is proved in \cite{FKRSp1} (and in  \cite[Section 2.2]{SaSaR}).
\begin{Pn} \label{PnWsead} The Weyl function always exists and it is unique. Moreover, the Weyl function
is contractive, that is, $\vp(z)^*\vp(z) \leq I_{m_1}$.
\end{Pn}
Inverse problems for the classical selfadjoint Dirac  systems had been actively studied since 1950s \cite{Kre1, LS}
and various interesting results were 
published last years (see, e.g., \cite{AHM, BrEKT, CG2, CGR, EGNT, GeGo, Hryn, MP, SaSaR}).
In particular, the inverse problem to recover the $m \times m$ locally square summable potential $V$ of  the Dirac system
\eqref{p1} from its Weyl function $\vp(z)$ was dealt with in \cite{Sa14}.
The case of rectangular (not necessarily square) matrix functions $v$ was studied there
and, moreover, only local square-summability of $v$ was required, which was essentially less
than in the preceding works.
Thus a problem formulated by F. Gesztesy was solved and interesting applications to the case
of Schr\"odinger-type operators with distributional matrix-valued potentials followed \cite{EGNST}.
\begin{Tm}\label{TmMai}\cite{Sa14}
Let Dirac system  \eqref{p1} be given on $[0, \, \infty)$,
let its  potential $V$ be locally square-summable and let
$\vp$ be the Weyl function of this system. Then $V$ is uniquely recovered
from $\vp$.
\end{Tm}
It is essential that a procedure to solve inverse problem is given in \cite{Sa14} (and will be formulated
in the next subsection). Therefore, we see that if we know the evolution $\vp(t,z)$ of the Weyl function
of the system $y_x=G(x,t,z)y$, where $G$ is given by \eqref{p5}, we may recover the solution $v(x,t)$
of  dNLS \eqref{p3}. The same scheme works for many other integrable equations.

\subsection{Recovery of the potential from the Weyl function}
In this subsection we describe a procedure 
to solve inverse problem for system \eqref{p1} with locally square summable potentials,
which is summed up in \cite[Theorem 4.4]{Sa14}.
The potential $v$ of  \eqref{p1} is easily expressed in terms of the block rows
$\b$ and $\g$ of the fundamental solution $u(x,z)$ (normalized by \eqref{p6}) at $z=0$:
\begin{align} &      \label{p8}
\b(x)=\begin{bmatrix}
I_{m_1} & 0
\end{bmatrix}u(x,0), \quad \g(x)=\begin{bmatrix}
0 &I_{m_2}
\end{bmatrix}u(x,0).
\end{align} 
Indeed, \eqref{p1} yields $u(x,z)^*ju(x, \ov{z})=j$. Hence, $u(x,0)ju(x, 0)^*=j$, and so
\begin{align} &      \label{p9}
\b j \b^*\equiv I_{m_1}, \quad \g j \g^*\equiv -I_{m_2}, \quad \b j \g^*\equiv 0. 
\end{align} 
Using \eqref{p1} and \eqref{p9} we derive
\begin{align} &      \label{p10}
v(x)=\I \b^{\prime}(x)j\g(x)^* \qquad \left(\b^{\prime}(x):=\left(\frac{d}{dx}\b\right)(x)\right).
\end{align}  
Moreover, relations \eqref{p9} and initial condition $\b(0)=\begin{bmatrix}I_{m_1} & 0\end{bmatrix}$
(which follows from \eqref{p6}, \eqref{p8}) provide a procedure to recover $\b$ from $\g$. Namely, we have
\begin{align} \label{p11}&
\b(x)= \b_1(x)\wt \b(x), \quad \wt \b(x):=\begin{bmatrix}
I_{m_1} & \g_1(x)^*(\g_2(x)^*)^{-1}
\end{bmatrix};
\\  \label{p12}&
\b_1^{\prime}=-\b_1(\wt \b^{\prime}j\wt \b^*)( \wt \b j \wt \b^*)^{-1}, \quad \b_1(0)=I_{m_1} ,
\end{align} 
where $\b_1$ is the left $m_1\times m_1$ block of $\b$ and $\g_1$ and $\g_2$ are the $m_2\times m_1$
and $m_2\times m_2$, respectively, blocks of $\g$. We note that \eqref{p9} yields $\det \g_2\not=0$. It remains to recover $\g$ from $\vp$. 
Now, we formulate the corresponding theorem.
\begin{Tm}\label{PrIP} Let Dirac system  \eqref{p1} be given on $[0, \, \infty)$,
let its  potential $V$ be locally square-summable and let
$\vp$ be the Weyl function of this system. 

Then $V$ is uniquely recovered
from $\g$ using formulas \eqref{p2} and \eqref{p10}, where $\b$ is given by \eqref{p11}
and $\b_1$ $($in \eqref{p11}$)$ is the unique solution of the first order linear differential system $($with initial condition$)$
in \eqref{p12}. In turn, the matrix function $\g$ is uniquely recovered from $\vp$ in three steps which are  described below.\\
$(i)$ Introduce a matrix function $\Phi_1(x)$   taking a Fourier transformation
\begin{align} \label{p13}&
\Phi_1\Big(\frac{x}{2}\Big)=\frac{1}{\pi}\E^{\eta x}{\mathrm{l.i.m.}}_{a \to \infty}
\int_{-a}^a\E^{-\I \xi x}\frac{\vp(\xi+\I \eta)}{2\I(\xi +\I \eta)}d\xi, \quad \eta >0,
\end{align} 
where ${\mathrm{l.i.m.}}$ stands for the entrywise limit in the norm of  $L^2(0,\ell)$, 
$0<\ell \leq~\infty$. Here $\Phi_1(x)$ $(0\leq x<\infty)$ is a well-defined differentiable matrix 
function $($with a locally square-summable derivative$)$, which does not depend on $\eta >0$. \\
$(ii)$ Introduce a family of pairs of operators $\Pi_{ l}\in B\Big(\BC^m,L^2_{m_2}(0,  l)\Big)$
and $S_{ l}$ acting in $L^2_{m_2}(0,  l)$  $\, (l\in \BR_+):$
\begin{align} \label{p14}&
\Pi_{ l}g=\Phi_1(x)g_1+g_2 \quad {\mathrm{for}} \quad g_k\in \BC^{m_k}, \, g= \begin{bmatrix} g_1\\ g_2 \end{bmatrix};
\\ \label{p15}&
S_{ l}=I-\frac{1}{2}\int_0^{ l}\int_{|x-t|}^{x+t}\Phi_1^{\prime}\left(\frac{ r+x-t}{2}\right)
\Phi_1^{\prime}\left(\frac{ r+t-x}{2}\right)^*d r \, \cdot \, dt.
\end{align} 
The defined above operator $S_{ l}$ belongs $B\Big(L^2_{m_2}(0,  l)\Big)$, it   is positive-definite $($i.e., $S_{ l}~>~0)$ and boundedly invertible,
and the matrix function $\Pi_{ l}^*S_{ l}^{-1}\Pi_{ l}$ is absolutely continuous with respect to $l$. \\
$(iii)$ The matrix function $($Hamiltonian$)$ $H( l)=\g( l)^*\g( l)$ is recovered via the formula
 \begin{align} &      \label{p16}
H( l)=\left(\Pi_{ l}^*S_{ l}^{-1}\Pi_{ l}\right)^{\prime}.
\end{align} 
From $H$ we recover first 
\begin{align} \label{p17}&
\g_2^{-1}\g_1=\left(\begin{bmatrix}
0 &I_{m_2}
\end{bmatrix}H\begin{bmatrix}
0 \\ I_{m_2}
\end{bmatrix}\right)^{-1}
\begin{bmatrix}
0 &I_{m_2}
\end{bmatrix}H\begin{bmatrix}
I_{m_1} \\ 0
\end{bmatrix}.
\end{align} 
Using $\g_2^{-1}\g_1$, we recover $\g_2$ from the equation
\begin{align} 
 &      \label{p18}
\g_2^{\prime}=
\g_2(\g_2^{-1}\g_1)^{\prime}(\g_2^{-1}\g_1)^*
\big(I_{m_2}-(\g_2^{-1}\g_1)(\g_2^{-1}\g_1)^*\big)^{-1}, \quad \g_2(0)=I_{m_2}.
\end{align} 
Finally, from $\g_2^{-1}\g_1$ and $\g_2$ the expression for $\g$
is immediate.
\end{Tm}
We note that indices in the notations  of operators (e.g., in the notations of operators $\Pi_{ l}$ and $S_{ l}$ above 
as well as in the notations $\, S_x$ and $\, S_l$ in Theorem \ref{wTmSk}) indicate the spaces, in which the operators act, and
should be distinguished from the indices meaning partial differentiation like in \eqref{p3}.
\subsection{Zero curvature representation as a compatibility \\ condition}\label{Comco}
In this subsection we do not require any special structure of $G$ and $F$ (like the structure in \eqref{p5})
and consider a general-type zero curvature representation \eqref{p4}.
The term "compatibility condition"
is very often used with respect to the representation \eqref{p4} and it is easy to derive \eqref{p4} 
from the existence of common fundamental solutions $w$ (i.e., from the compatibility) of systems
 \begin{align}&      \label{p19}
w_x=Gw, \quad w_t=Fw.
\end{align}
In other words, it is easy to show that \eqref{p4}
 is a necessary compatibility condition for systems \eqref{p19}.
However, the proof  of sufficiency is more complicated. It appeared first in \cite{SaL20, SaLev2}
and in greater detail and generality in \cite{SaAF, SaSaR}.

More precisely, assuming  that \eqref{p4} holds  in the semistrip
  \begin{align}&
   \label{p20}
\Om_a=\{(x,\, t):\,0 \leq x <\infty, \,\,0\leq t<a\}
\end{align}
we introduce solutions of systems $W_x=GW$ with fixed values of $t$ and of systems $R_t=FR$ with
fixed values of $x$:
\begin{align}& \label{p21}
W_x(x,t,z)=G(x,t,z) W(x,t,z), \quad W(0,t,z)=I_m;
\\ & \label{p22}
R_t(x,t,z)=F(x,t,z) R(x,t,z), \quad R(x,0,z)=I_m.
\end{align}
Next, we formulate \cite[Theorem 6.1]{SaSaR}.
\begin{Tm}\label{TmComp}
Let $m \times m$ matrix functions $G$ and $F$ and their derivatives $G_t$ and
$F_x$  exist on the semistrip $\Om_a$,
let $G$, $G_t$, and $F$ be continuous with respect to $x$ and $t$ on $\Om_a$,
 and let \eqref{p4} hold. Then we have the equality 
\begin{equation} \label{p23}
W(x,t,z)R(0,t,z)=R(x,t,z)W(x,0,z).
\end{equation}
\end{Tm}
Thus, \eqref{p4} implies \eqref{p23}, and so
 \begin{align}&
  \label{pd}  
w(x,t,z):=W(x,t,z)R(0,t,z)=R(x,t,z)W(x,0,z)
\end{align}
is the fundamental solution
of systems \eqref{p19} normalized by the initial condition 
 \begin{align}&
  \label{p24}  
 w(0,0,z)=I_m.
\end{align}

The problem of  compatibility  of the systems \eqref{p19} may be formulated
as an initial-boundary value problem. Indeed, normalizing $w$ in \eqref{p19} via \eqref{p24}
and taking into account \eqref{p21} and \eqref{p22} we obtain initial-boundary conditions
on $w$:
  \begin{align}&
  \label{p25}  
 w(x, 0,z)=W(x, 0,z), \quad w(0,t,z)=R(0,t,z).
\end{align}
From this point of view, formula \eqref{pd} gives
a solution of the initial-boundary value problem \eqref{p19}, \eqref{p25}.
This solution plays a crucial role in the study of the evolution of Weyl functions.
\subsection{Evolution of Weyl functions}\label{Evol} 
Inverse Spectral Transform (instead of  the Inverse Scattering Transform) was first
used for initial-boundary value problems in \cite{Ber, KvM}. More precisely,  a special kind of initial-boundary value problem for Toda lattice
with a linear law of evolution of the spectral and Weyl functions was studied in \cite{Ber, KvM}. In the papers
\cite{SaL20, SaLev2}, an essentially more general case of the initial-boundary value problem
for Toda lattice as well as some initial-boundary value problems for continuous integrable
systems (including square matrix dNLS) were dealt with, and the law of evolution of the Weyl
function was presented in the form of M\"obius transformation. See further
results and references, for instance, in \cite{SaA021, SaAg, SaSaR}. 

In this section we 
demonstrate Inverse Spectral Transform   approach on the important case of dNLS \eqref{p3}, where
$v$ is an $m_1\times m_2$ matrix function. For that purpose,  we use $W$ and $R$ given by \eqref{p21}
and \eqref{p22} assuming  that $G$ and $F$ have the form \eqref{p5}. First, we omit $t$ in the notations
and introduce the set (Weyl circle) $\cln(b,z)$ of functions of the form
\begin{equation} \label{p26}
\phi(b,z,\clp)=\begin{bmatrix}
0 &I_{m_2}
\end{bmatrix}  W(b,z)^{-1}\clp(z)
 \Big(\begin{bmatrix}
I_{m_1} & 0
\end{bmatrix}W(b,z)^{-1}\clp(z)\Big)^{-1},
\end{equation}
where $\clp(z)$ are $m \times m_1$ nonsingular  meromorphic matrix functions 
 with property-$j$, that is,
 \begin{align}\label{p27}&
\clp(z)^*\clp(z) >0, \quad \clp(z)^* j \clp(z) \geq 0 \quad (z\in \BC_+).
\end{align}
The next formula gives an important property of Weyl functions $\vp(z)$:
 \begin{align}&      \label{p28}
\vp(z)=\lim_{b \to \infty}\phi(b,z)
\end{align} 
for any set of functions $\phi(b,z)\in \cln(b,z)$. In order to derive an expression for the Weyl
function $\vp(t,z)$ of the system $y_x=G(x,t,z)y$, we rewrite \eqref{p23} in the form of the equality
$$W(b,t,z)^{-1}=R(0,t,z)W(b,0,z)^{-1}R(b,t,z)^{-1}$$
and substitute this equality into \eqref{p26}.  Now, passing  to the limit 
and using \eqref{p28} we obtain the required expression for $\vp(t,z)$. That is, we obtain the dependence of the
Weyl functions on the parameter $t$ (i.e., evolution of the Weyl function) in the case of dNLS.
\begin{Tm}\label{TmevNLS} \cite{Sa14b}
Let an $m_1 \times m_2$ matrix function $v(x,t)$ 
 be  continuously differentiable on the semistrip $\Om_a$
 and let  $v_{xx}$ exist. 
 Assume that $v$  satisfies the dNLS equation \eqref{p3}
   as well as  the following inequalities
 $($for all $0 \leq t<a$ and some values $M(t)\in \BR_+):$
  \begin{equation} \label{p29}
  \sup_{x \in \BR_+, \, 0\leq s \leq t}\|v(x,s)\| \leq M(t).
\end{equation}
 Then the evolution $\vp(t,z)$ of the Weyl functions  of  Dirac systems \\  $y_x(x,t,z)=G(x,t,z)y(x,t,z)$ is given $($for $z\in \BC_+)$ by the equality
 \begin{equation} \label{p30} 
\vp(t,z)= \big(R_{21}(t,z)+R_{22}(t,z)\vp(0,z)\big)
\big(R_{11}(t,z)+R_{12}(t,z)\vp(0,z)\big)^{-1},
\end{equation}
where $R_{ik}(t,z)$ are the $m_i\times m_k$ blocks of $R(0,t,z)$.
\end{Tm}

\section{Generalized Weyl functions} \label{GW}
\setcounter{equation}{0} 
\subsection{Skew-selfadjoint Dirac
system} \label{Sk}
 The system
\begin{equation}       \label{s1}
\frac{d}{dx}y(x, z )=(\I z j+jV(x))y(x,
z ) \qquad
(x \geq 0, \quad z \in \BC),
\end{equation}
where $j$ and $V$ have the same form \eqref{p2} as in the Dirac system \eqref{p1},
is called skew-selfadjoint Dirac.
Weyl theory of skew-selfadjoint Dirac systems was studied in the papers \cite{CG2, FKRSsk, SaA021, SaAg}.
Some further references as well as the results of this subsection are contained in  \cite[Ch.3]{SaSaR}.
Like in the case of selfadjoint Dirac system, the fundamental solution of system
\eqref{s1} is denoted by $u(x,z)$, and this solution is normalized by the condition \eqref{p6}.
The notation $\BC_M$ stands for the open half-plane $\{z:\, \Im (z)>M>0\}$.
\begin{Dn} \label{DnW2} An $m_2 \times m_1$ matrix function $\vp$, which is  holomorphic in $\BC_M$
$($for some  $M>0)$ and satisfies the inequality
\begin{align}&      \label{s2}
\int_0^{\infty}
\begin{bmatrix}
I_{m_1} & \vp(z)^*
\end{bmatrix}
u(x,z)^*u(x,z)
\begin{bmatrix}
I_{m_1} \\ \vp(z)
\end{bmatrix}dx< \infty , \quad z\in \BC_M,
\end{align}
is called a Weyl function of the skew-selfadjoint Dirac system  \eqref{s1}.
\end{Dn}
First, assume that $v$ is bounded, that is,
 \begin{align}\label{s3}
 \|v(x)\| \leq M  \quad \mathrm{for} \,\, x\in [0, \, \infty).
 \end{align}
 For $z\in \BC_M$ (with $M$ in $\BC_M$ given in \eqref{s3}), we have the following proposition.
\begin{Pn}\label{PnWgW} Let Dirac system  \eqref{s1} be given on $[0, \, \infty)$ and assume that \eqref{s3} holds. Then there is a unique
Weyl function $\vp(z)$ $($for $z\in \BC_M)$ of this system.  Moreover, $\vp(z)$ is contractive in $\BC_M$ and the inequality
\begin{align}&      \label{s4}
\sup_{x \leq l, \, z\in \BC_M}\left\| \E^{-\I z x}u(x,z)\begin{bmatrix}
I_{m_1} \\ \vp(z)
\end{bmatrix} \right\|<\infty
\end{align}
holds on  any finite interval $[0, \, l]$.
\end{Pn}
\begin{Dn}   \label{DnGW}
 A   generalized Weyl function $(${\rm GW}-function$)$ of  the system
\eqref{s1}, where $v$ is locally bounded on $[0, \, \infty)$,
is an  $m_2 \times m_1$ matrix function $\vp$ such that for some $M>0$ it is analytic in $\BC_M$ and the
inequality \eqref{s4}  holds for each
$l < \infty$.
\end{Dn}
Our next proposition (together with Proposition \ref{PnWgW}) shows that, in the case
of $v$ bounded on the semiaxis, the definitions of Weyl function and GW-function are equivalent,
that is, GW-function coincides with Weyl function.
 \begin{Pn}  \label{PnUGW} For any system
\eqref{s1}, where $v$ is locally bounded on $[0, \, \infty)$, there is no more than one {\rm GW}-function.
  \end{Pn}
The following theorem (see \cite[Theorem 3.30]{SaSaR}) is necessary in the process of solving Goursat problem for the sine-Gordon equation in Section \ref{SG}.
The procedure to construct $\M(\vp)$ (solve inverse problem) given in this theorem is in many respects similar to the corresponding procedure for the selfadjoint Dirac system but
we do not require here apriori that $\vp$  is a GW-function. We note that $\b$ and $\g$ stand again for  the block rows of $u$,  that is, \eqref{p8} holds.
 \begin{Tm}\label{wTmSk} Let an $m_2 \times m_1$ matrix function $\vp(z)$ be holomorphic in $\BC_M $ $($for some $M>0)$ and satisfy condition
 \begin{align}\label{s5}
 \sup\|z^2(\vp(z)-\phi_0/z)\|< \infty \quad (z \in \BC_M),
 \end{align}
 where $\phi_0$ is an $m_2 \times m_1$ matrix. Then $\vp$ is a GW-function
 of a skew--selfadjoint Dirac system.
 This Dirac system is uniquely recovered from $\vp$ using the following procedure.\\

$(i)$ First recover the matrix function $\Phi_1:$ 
 \begin{align} \label{s6}&
\Phi_1\Big(\frac{x}{2}\Big)=\frac{1}{\pi}\E^{\eta x}{\mathrm{l.i.m.}}_{a \to \infty}
\int_{-a}^a\E^{-\I \xi x}\frac{\vp(\xi+\I \eta)}{2\I(\xi +\I \eta)}d\xi, \quad \eta >M,
\end{align}
where {\rm l.i.m.} stands for the entrywise limit in the norm   $L^2(0,\ell)$ $(0<\ell \leq~\infty)$.

$(ii)$ Next, for the values $l\in \BR_+$, introduce operators $S_l \in B\left(L^2(0,l)\right):$
\begin{align}  &
\label{s7}
S_l= I +\int_0^ls(x,t)\, \cdot \, dt, \quad s(x,t)=\int_0^{\min(x,t)}\Phi_1^{\prime}(x-r)\Phi_1^{\prime}(t-r)^*dr.
\end{align}
These operators are well-defined, positive definite and boundedly invertible.

$(iii)$ We recover  $\b(x)$ via the formula $:$
 \begin{align}&      \label{s8}
\b(x)=\begin{bmatrix}I_{m_1} &0 \end{bmatrix}
-\int_0^x\Big(S_x^{-1}\Phi_1^{\prime}\Big)(t)^*\begin{bmatrix}\Phi_1(t) & I_{m_2} \end{bmatrix}dt,
\end{align} 
where $S_x^{-1}$ is applied to $\Phi_1^{\prime}$ columnwise. Since $\b \b^* \equiv I_{m_1}$, 
we can construct a  differentiable matrix function $\wt \g$  $($with locally bounded derivative$)$ such that
\begin{align}     \label{s9}&
\b \wt \g^*\equiv 0, \quad \wt \g \wt \g^*>0, \quad \wt \g(0)=\begin{bmatrix}
0 & I_{m_2}\end{bmatrix}.
\end{align}
Then $\g=\wt \vt \wt \g$, where $\wt \vt$ is determined by the equation and initial condition below $:$
\begin{align}     \label{s10}&
\wt \vt^{\prime}=-\wt \vt\wt \g^{\prime} \wt \g^* (\wt \g \wt \g^*)^{-1}, \quad \wt\vt(0)=I_{m_2}.
\end{align}

$(iiii)$ Finally, $v$ is given by
\begin{align}     \label{s11}&
v(x)=\b^{\prime}(x)\g(x)^*,
\end{align}
$V$ has the form \eqref{p2}, and both $v$ and $V$ are locally bounded.
\end{Tm}
\subsection{Linear system auxiliary to the $N$-wave equation} \label{Nw}
Nonlinear optics ($N$-wave) equation 
\begin{align}  \label{NW1} &
\big[D,  \vr_t \big] - \big[\wh D,  \vr_x \big]
=
[D, \vr]\, [\wh D, \vr] - [\wh D, \vr]\, [ D, \vr], \
\end{align}
where $\vr(x,t)=\vr(x,t)^*$ is an $m\times m$ matrix function, $[D, \vr]=D\vr - \vr D$ and
\begin{align} 
\label{NW2} &  D= 
\di  \{  d_1,   d_2, \ldots,
  d_m \}, \quad d_1>d_2>\ldots>d_m>0, \\
\label{NW3} &   
  \wh D=\wh D^*=\di\{\wh d_1, \wh d_2, \ldots, \wh d_m\},
\end{align}
admits zero curvature representation \eqref{p4}, where $G$ and $F$ have the form
\begin{align}& \label{NW4}
G(x,t,z)=\I z D-\zeta(x,t), \quad F(x,t,z)=\I z \wh D-\wh \zeta(x,t); \\ 
& \label{NW4'}
\zeta:=[D, \vr], \quad \wh \zeta: =[\wh D, \vr];
\end{align}
see \cite{ZaMa} for the case $N=3$ and \cite{AbH} for $N>3$.

We shall need
some preliminary results on the Weyl theory of the auxiliary system $y_x=Gy$ $(x \geq 0)$
from \cite[Ch. 4]{SaSaR}, see also \cite{SaARMS, SaA17}. The normalized fundamental solution $u$
of such system is defined by the formula
\begin{align}& \label{NW5}
\frac{d}{d x}u(x,z)=\big(\I z D-\zeta(x)\big)u(x,z), \quad u(0,z)=I_m \quad (\zeta =- \zeta^*).
\end{align}
Here and further we assume that $D$ is a fixed matrix satisfying \eqref{NW2}. 
\begin{Dn}   \label{NWDn2}
 A  generalized Weyl function $($GW-function$)$ of system
\eqref{NW5},  where $\zeta$ is locally bounded, is an  $m \times m$ matrix function $\vp$ such that for some $M>0$ it is analytic in  the domain
$\BC^{-}_M = \{ z \, \colon \, \Im(z) < - M \}$ and the
inequality  
\begin{equation}
      \sup\limits_{x \, \leq \, l, \ \Im (z) <  -M} \,
          \bigl\|
                u(x,z)
                \varphi (z)\exp  \{-\I z x D\}
          \bigr\|
    < \infty 
   \label{NW6}
\end{equation}
holds for each
$l < \infty$.
\end{Dn}  
The inverse spectral problem  for  system
\eqref{NW5} is the problem to recover (from an
analytic matrix function $\varphi(z)$)
a locally bounded potential
$\zeta (x) = - \zeta(x)^*$
$(\zeta_{kk} \equiv 0)$ such that  $\vp$ is a $GW$-function of the corresponding
system \eqref{NW5}, that is,
\eqref{NW6}  is valid. The notation $\mf M$
stands for the solution of this inverse problem, that is, for
an operator mapping the pair $D$ and
$\varphi $
into the corresponding potential
$\zeta$
$($i.e., ${\mf M} (D, \varphi) = \zeta)$.   
\begin{Tm}\label{NWTmUniq}
For any
matrix function
$\varphi (z)$ which is
analytic and bounded in
$\BC^{-}_M = \{ z \, \colon \, \Im(z) < - M \}$
and has the
property
\begin{equation}
      \int\limits_{- \infty}^\infty
          \bigl(
               \varphi (\xi+\I \eta) - I_m^{\,}
          \bigr)^* \,
          \bigl(
               \varphi (\xi+\I \eta) - I_m^{\,}
          \bigr)
      \, d \xi
      < \infty
\qquad
      (
       \eta <- M),
    \label{NW8!}
\end{equation}
there is at most one solution of the inverse spectral problem.
\end{Tm}
An existence condition (for the solution of the inverse spectral problem) is given in the next theorem. 
\begin{Tm}\label{NWInvPr} Let $\vp(z)$ be analytic in $\BC^{-}_M $ for some $M>0$
and satisfy in $\BC^{-}_M$ the inequalities
\begin{align} & \label{NW9}
\sup_{z\in \BC^{-}_M} 
            \big\| \,
           z ( \varphi (z) - I_m) \,
           \big\|
            < \infty, \quad \det \, \varphi (z) \not= 0.
\end{align}
Assume also that for   some matrix $\phi_0$ and for all fixed $\eta<-M$ we have 
\begin{align} & \label{NW10}
           (\xi+\I\eta)\big( \varphi (\xi+\I\eta) - I_m - \phi_0 / (\xi+\I \eta) \big)
            \in L_{m \times m}^2 (- \infty, \infty).
  \end{align}
Then the solution of the inverse problem exists and  is unique. 
\end{Tm}
\begin{Rk} \label{RkPrM}
The procedure to construct ${\mf M} (D, \varphi)$ under conditions of Theorem \ref{NWInvPr} is given in \cite[Theorem 4.10]{SaSaR}.
This procedure is similar to the procedures in Theorems \ref{PrIP} and \ref{wTmSk}.
\end{Rk}

Further in this paper, we deal with the case of $\zeta$ bounded on $[0,\infty)$:
\begin{equation}
      \sup\limits_{0 < x < \infty} 
          \bigl\|
                \zeta (x)
          \bigr\|
      < \infty
    \label{NW7},
\end{equation}
 and  $GW$-functions $\vp(z)=\{\vp_{ij}(z)\}_{i,j=1}^m$ normalized by
\begin{align}&
         \varphi_{ij}^{\,} (z) \equiv 1 \quad  \mbox{   for   } \quad   i= j, \qquad
         \varphi_{ij}^{\,} (z) \equiv 0 \quad \mbox{   for   } \quad  i > j.
    \label{NW8}
\end{align} 
The next proposition follows from \cite[Subsections 4.1.1 and 4.1.3]{SaSaR}.

\begin{Pn} When \eqref{NW7} holds, a normalized $($by \eqref{NW8}$)$ $GW$-function $\vp(z)$ of \eqref{NW5}
exists and is unique.
\end{Pn}

\section{Sine-Gordon theory in a semistrip} \label{SG}
\setcounter{equation}{0} 
{\bf 1.} Sine--Gordon equation  in the light cone coordinates (SGE) has the form
\begin{equation}      \label{SG1}
\frac{\prt^2}{\prt t \prt x}\psi=2\sin (2\psi).
\end{equation}
Local solutions of the Goursat problem for SGE were studied in \cite{Kri, LeSa}
and global solutions were constructed first in \cite{SaAg}. The results of this section
are obtained mostly in \cite{SaAg} with some modifications and developments in
\cite[Section 6.2]{SaSaR} (see also further references in \cite{SaSaR}).
SGE admits zero curvature representation \eqref{p4}, where $G$ and $F$ have the form
 \cite{AKNS0}:
\begin{align}      \label{SG2}
& G(x, t,z)=\I zj+jV(x,t), \quad V=\left[
\begin{array}{lr}
0& v \\   v & 0
\end{array}
\right], \quad v=-\frac{\prt}{\prt x}\psi,
\quad \psi=\ov \psi,
\\
& \label{SG3}
 F(x, t,z)=\frac{1}{\I z}
\left[
\begin{array}{lr}
\cos (2\psi)&  \sin (2 \psi) \\ \sin (2 \psi) &-\cos (2 \psi)
\end{array}
\right],
\end{align}
and $j$ has the form \eqref{p2} with $m_1=m_2=1$. Hence, system $y_x=Gy$ is a skew-selfadjoint Dirac system \eqref{s1}, where $m_1=m_2=1$. 
The evolution $\vp(t,z)$ of its $GW$-function
is given by the formula
\begin{equation} \label{SG4}
\vp(t,z)=\frac{R_{21}(t,z)+R_{22}(t,z)\vp(0,z)}
{R_{11}(t,z)+R_{12}(t,z)\vp(0,z)},
\end{equation}
where, as usual, $R_{ik}(t,z)$ are the blocks of $R(0,t,z)$ determined by \eqref{p22}.
If  the boundary value
$\psi(0,t)$ is given,
using \eqref{p22} and \eqref{SG3} we recover $R_{ik}(t,z)$. Let us formulate the evolution result rigorously.
\begin{Tm}\label{SGevolSG}
Let $\psi(0,t)$ and $ \psi_x(x,t)$ be continuous functions on $[0, \, a)$
and
$\Om_a$, respectively, and let $ \psi_{xt}(x,t)$ exist.  
Let also SGE \eqref{SG1} hold on $\Om_a$, and assume that 
 $\vp(0,z)$ is the $GW$-function
 of the system $ y_x=G(x,0,z)y$, where $G$ is given by \eqref{SG2}.

 Then, the function $\vp(t,z)$ of the form  \eqref{SG4} is the $GW$-function  of the system $ y_x=G(x,t,z)y$.
\end{Tm}
If the initial--boundary  conditions
\begin{equation} \label{SG5}
\psi(x,0)=h_1(x), \quad \psi(0,t)=h_2(t),  \quad h_1(0)=h_2(0)
\end{equation}
($h_k=\ov{h_k}$ for $(k=1,2)$), are given, we  recover $\vp(0,z)$ from the first
condition, recover $R_{ik}(t,z)$ from the second  condition, and construct $\vp(t,z)$ using \eqref{SG4}. 
Then we recover $v(x,t)$ using the procedure
to construct the solution $\M(\vp)$ of the inverse problem, which is described in Theorem \ref{wTmSk}.
From $v$ we immediately recover $\psi$ and show that $\psi$ satisfies \eqref{SG1} and \eqref{SG5}
(see the proof of \cite[Theorem 6.19]{SaSaR}). Thus, an existence theorem follows.
\begin{Tm}\label{SGWsg}  
Assume that  $h_1$ is boundedly  differentiable on all the finite intervals on $[0,\, \infty)$ and that
$h_2$ is continuous on $[0, \, a)$.

Moreover, assume that the $GW$-function $\vp(0,z)$ of the system \eqref{s1}, where
\begin{equation} \label{SG6}
m_1=m_2=1, \quad V(x)=-\left[
\begin{array}{lr}
0& h_1^{\prime}(x) \\  h_1^{\prime}(x) & 0
\end{array}
\right],
\end{equation}
exists
 and satisfies  \eqref{s5}. Then a solution of the
initial-boundary value problem \eqref{SG1}, \eqref{SG5}  exists
and is given by the equality
\begin{equation} \label{SG7}
\psi(x,t)=h_2(t)-\int_0^{x}\Big(\M\big(\vp(t,z)\big)\Big)(\xi)d\xi,
\end{equation}
where $\vp(t,z)$ is obtained from \eqref{SG4}.
\end{Tm}
\begin{Rk}
Under conditions of  Theorem \ref{SGWsg}, for each $0<c<a$ there is $M(c) >0$ such that all functions $\vp(t,z)$ $(0\leq t \leq c)$      satisfy $GW$-function
requirement \eqref{s4} and asymptotic inequality \eqref{s5} in the same half-plane $\BC_{M(c)}$.
\end{Rk}

Some sufficient conditions on $h_1$, under which the requirements on $\vp(0,z)$ in Theorem \ref{SGWsg}
hold, were derived using important paper \cite{BC0} (see \cite[Corollary~6.21]{SaSaR}).  We formulate these
conditions.
\begin{Cy} \label{SGMCy}  If  $v \in L^1(\BR_+)$, then there is a $GW$-function $\vp$
of  the system \eqref{s1} $($where $V$ has the form \eqref{p2} and $m_1=m_2=1)$. 
Moreover, if  $v$ is two times differentiable and $v, \, v^{\prime}, \, v^{\prime \prime}
\in L^1(\BR_+)$, then this GW-function $\vp$ satisfies the asymptotic condition \eqref{s5}.
\end{Cy}
{\bf 2.} Complex sine-Gordon equation
was introduced (and its integrability was treated) \cite{LuRe, Po} only several
years after the seminal paper \cite{AKNS0} on the integrability of SGE was published.  Complex sine-Gordon equation is  more general than SGE
and has the form
\begin{equation}      \label{SG8}
\psi_{xt}+4 \cos \psi (\sin \psi)^{-3} \chi_x \chi_t =2\sin (2\psi), 
\quad  \chi_{xt}=(2/\sin (2 \psi)) (\psi_x \chi_t+\psi_t \chi_x),
\end{equation}
where $\psi=\ov\psi$ and $\chi=\ov\chi$.
There are also two constraint equations
\begin{equation}      \label{SG9}
2( \cos \psi)^2\chi_x-( \sin \psi)^2 \om_x=2c( \sin \psi)^2, \quad
2( \cos \psi)^2 \chi_t+( \sin \psi)^2 \om_t=0, 
\end{equation}
where $c$ is a constant $( c=\ov c \equiv {\mathrm{const}})$
and $\om=\ov \om$.
For further developments and applications of the results on \eqref{SG8}, \eqref{SG9} see, for instance,
some discussions and references in   \cite{BP, BoTz, SaSaR}.  Below we formulate several results from
\cite[Subsection 6.2.1]{SaSaR}.

If $\sin(2\psi)\not=0$ and \eqref{SG8} and \eqref{SG9} hold, then zero curvature equation \eqref{p4}, where
\begin{align}   \label{SG10}
&G(x,t,z):=\I zj+jV(x,t), \quad F(x,t,z):=-\frac{\I}{z+c}\vt(x,t)^*j\vt(x,t),
\\  \label{SG11}&
 j =\left[
\begin{array}{cc}
1& 0 \\ 0 & -1
\end{array}
\right], \quad V=\left[
\begin{array}{lr}
0& v \\  \ov v & 0
\end{array}
\right], \quad v=\left(-\I \frac{\prt \psi}{\prt x}+2(\cot \psi) \frac{\prt \chi}{\prt x}\right)\E^{\I(\om-2\chi)},
\\ & \label{SG12}
\vt(x,t): = D_1(x,t)\left[\begin{array}{cc}
\cos\psi& \I \sin \psi \\
 \I \sin \psi &\cos\psi
\end{array}\right]D_2(x,t),\\
& \label{SG13}
D_1=\exp\big\{\I\big(\chi +({\om}/{2})\big)j\big\}, \quad D_2=\exp\big\{\I\big(\chi -({\om}/{2})\big)j\big\},
\end{align}
also holds. Now, we introduce initial-boundary conditions:
\begin{equation}\label{SG14}
v(x,0)=h_1(x), \quad \psi(0,t)=h_2(t), \quad \chi(0,t)=h_3(t), \quad \om(0,0)=h_4.
\end{equation}
Evolution of the $GW$-function $\vp$ is given in the next theorem.
\begin{Tm} \label{SGevol} Let  $\{\psi(x,t), \, \chi(x,t), \, \om(x,t)\}$ be a triple of real-valued 
and twice continuously differentiable functions on $\Om_a$.
 Assume that $\sin (2\psi) \not=0$, that $v$ given by \eqref{SG11} is bounded on $\Om_a$,
 and that complex sine-Gordon  \eqref{SG8}, \eqref{SG9} and conditions \eqref{SG14} hold.
 
Then the $GW$-functions $\vp(t,z)$ of the auxiliary skew-selfadjoint Dirac systems
$y_x=Gy$, where $G(x,t,z)$  is defined via \eqref{SG10} and \eqref{SG11},
exist in some $\BC_M$ and have the form \eqref{SG4}.

Here $R(t,z)=R(0,t,z)=\big\{R_{ik}(t,z) \big\}_{i,k=1}^2$ is defined by the equalities
\begin{align}\nn &
\frac{d}{d t}R(t,z)= \frac{1}{\I (z+c)}\E^{-\I d(t)j}
\left[
\begin{array}{lr}
\cos (2h_2(t))& \I \sin (2 h_2(t)) \\ -\I \sin (2 h_2(t)) &-\cos (2 h_2(t))
\end{array}
\right] 
\E^{\I d(t)j} R(t,z), 
\\
 \nn &
 R(0,z)=I_2, \quad d(t):=h_3(0)-\frac{1}{2}h_4+\int_0^t h_3^{\prime}(\xi)\big(\sin h_2(\xi)\big)^{-2}d\xi;
\end{align}
and $\vp(0,z)$ is the $GW$-function of  \eqref{s1} where $m_1=m_2=1$ and $v=h_1$.
\end{Tm}
\begin{Cy}\label{SGUniq} There is at most one triple $\, \{\psi(x,t), \, \chi(x,t), \, \om(x,t)\} \,$  of real-valued 
and twice continuously differentiable functions on $\quad \Om_a \quad$ such that $ \sin (2\psi)~\not=~0$, that $v$ is bounded, and
that complex sine-Gordon \eqref{SG8}, constraints \eqref{SG9} and initial-boundary conditions \eqref{SG14} are satisfied.
\end{Cy}

\section{Reduction of the initial-boundary conditions  and unbounded solutions} \label{Reduc}
\setcounter{equation}{0} 
\subsection{Reduction of the initial-boundary conditions in a quarter-plane}
{\bf 1.} In this subsection we consider integrable wave equations in the quarter-plane
  \begin{align}&
   \label{r1}
 \Om_{\infty}=\{(x,\, t):\,0 \leq x <\infty, \,\,0\leq t<\infty\}.
\end{align}
First, we consider the defocusing NLS  \eqref{p3} with the initial-boundary conditions
  \begin{align}&
   \label{r2}
v(x,0)=h_1(x), \quad v(0,t)=h_2(t), \quad v_x(0,t)=h_3(t).
\end{align}
Given $h_1(x)$, we recover $\vp(0,z)$ via \eqref{p28} and given
$h_2(t)$ and $h_3(t)$ we recover $R(0,t,z)$ (i.e., $R_{ik}(t,z)$ where $i,k=1,2$)
using \eqref{p22}. Hence, given the initial-boundary conditions $h_k$ ($k=1,2,3$)
we recover the right-hand side of \eqref{p30} and obtain
the evolution $\vp(t,z)$ of the Weyl function.
\begin{Tm}\label{TmQP}\cite{Sa14b} Let $v$  satisfy the conditions of Theorem    \ref{TmevNLS} 
and boundary conditions given by the second and third equalities in \eqref{r2}. Moreover, let
the boundary value functions $h_2$ and $h_3$ be continuous and bounded, i.e.,
  \begin{align}&
   \label{r3}
\sup_{0\leq t<\infty}\|h_2(t)\|<\wh M, \quad \sup_{0\leq t<\infty}\|h_3(t)\|<\breve M
\end{align}
for some $\wh M, \breve M \in \BR_+$. Then, in the domain
\begin{align}& \label{r4}
\cld=\{z: \quad \Im(z) \geq 1/2, \quad \Re(z)\leq - \wh M\}
\end{align}
we have the equality
\begin{align}& \label{r5}
\vp(0,z) = -\lim_{t \to \infty} R_{22}(t,z)^{-1}R_{21}(t,z).
\end{align}
\end{Tm}
{\bf Scheme of the proof.}  The contractiveness of $\vp(t,z)$ yields
\begin{align}& \label{r6}
\begin{bmatrix} I_{m_1} & \vp(t,z)^*  \end{bmatrix}j\begin{bmatrix} I_{m_1} \\ \vp(t,z)  \end{bmatrix} \geq 0.
\end{align}
Putting $R(t,z):=R(0,t,z)$ and using formula \eqref{p30} for $\vp(t,z)$, we derive 
\begin{align}& \label{r7}
 R(t,z)\begin{bmatrix} I_{m_1} \\ \vp(0,z)  \end{bmatrix}=\begin{bmatrix} I_{m_1} \\ \vp(t,z)  \end{bmatrix}(R_{11}(t,z)+R_{12}(t,z)\vp(0,z)).
\end{align}
It is immediate from \eqref{r6} and \eqref{r7} that
\begin{align}& \label{r8}
\begin{bmatrix} I_{m_1} & \vp(0,z)^*  \end{bmatrix}R(t,z)^* jR(t,z)\begin{bmatrix} I_{m_1} \\ \vp(0,z)  \end{bmatrix}\geq 0.
\end{align}
Moreover, the first inequality in \eqref{r3} implies that
\begin{align}& \label{r9}
\frac{d}{dt}\big(-R(t,z)^*jR(t,z)\big) \geq \wh M R(t,z)^*R(t,z), \quad z\in \cld.
\end{align}
From \eqref{r8} and \eqref{r9} we see that
\begin{align} & \label{r10}
\int_0^{\infty}\begin{bmatrix} I_{m_1} & \vp(0,z)^*
\end{bmatrix}R(s,z)^*R(s,z)
\begin{bmatrix} I_{m_1} \\ \vp(0,z)
\end{bmatrix}ds
\leq (1/\wh M) I_{m_1}.
\end{align}
After some considerations, using the boundedness of the left-hand side of \eqref{r10} and the inequalities
\eqref{r3}, one can show that
\begin{align}& \label{r11}
\lim_{t\to \infty}\left\| R(t,z)
\begin{bmatrix} I_{m_1} \\ \vp(0,z)
\end{bmatrix}
\right\|=0.
\end{align}
On the other hand, it easily follows from \eqref{r9} that 
$$R_{22}(t,z)^*R_{22}(t,z)\geq I_{m_2}.$$ 
Therefore,
formula \eqref{r11} implies \eqref{r5}.  An analogue of \eqref{r5} for a scalar dNLS
(and with a somewhat different proof) appeared already
in \cite{SaL20}.

In view of analyticity of Weyl functions, formula \eqref{r5}  means that $\vp(0,z)$, $\vp(t,z)$,
and thus also $v(x,t)$, may be recovered from the boundary conditions. \\
{\bf 2.} A similar result is valid for the focusing nonlinear Schr\"odinger equation:
\begin{align}\label{fNLS}
2v_t+\I(v_{xx}+2vv^*v)=0 ,
\end{align}
where $v$ is an $m_1\times m_2$ matrix function. Equation \eqref{fNLS} admits \cite{ZS1,ZS2} representation \eqref{p4}
where
\begin{align}  \label{r15} &
G=\I zj+jV, \quad
F=\I\big(z^2 j-\I z jV-\big(V_x+jV^2\big)/2\big),
\end{align}
and $j$ and $V$ are defined in \eqref{p2}. Evolution of the Weyl function is described in this case
by the following theorem.
  \begin{Tm}\label{evol} \cite{FKRSsk} Let an $m_1 \times m_2$ matrix function $v(x,t)$ 
 be  continuously differentiable on $\Om_a$ 
 and let  $v_{xx}$ exist. 
 Assume that $v$  satisfies the  focusing nonlinear Schr\"odinger equation \eqref{fNLS}  as well as  the following inequalities $:$
  \begin{equation} \nn
  \sup_{(x,t)\in {\Om_a}}\|v(x,t)\| \leq M, \quad \sup_{(x,t)\in {\Om_c}}\|v_x(x,t)\|<\infty
  \quad \mathrm{for}\,\mathrm{each} \quad 0< c<a.
\end{equation}
 Then, the evolution $\vp(t,z)$ of the $GW$-functions  of   the skew-selfadjoint Dirac systems $y_x=Gy$
  is given $($for $z \in \BC_M)$ via the M\"obius transformation \eqref{p30},
 where the coefficients $R_{ik}$ are determined by the formula \eqref{p22} and the equality $\{R_{ik}(t,z)\}_{i,k=1}^2:=R(0,t,z)$. Here $G$ and $F$ have the form \eqref{r15}.
  \end{Tm}
\begin{Rk}
It is immediate that if the conditions of Theorem \ref{evol} hold for each $a<\infty$, then
\eqref{p4} describes the evolution of $\vp(t,z)$ for all $0\leq t <\infty$.
\end{Rk}
In order to derive an analogue of Theorem \ref{TmQP} we shall require boundedness of $v$ on $\Om_{\infty}$ (instead of boundedness
on all $\Om_a$) and boundedness of $v_x(0,t)$ on $\BR_+$:
\begin{align}& \label{r16}
  \sup_{(x,t)\in {\Om_{\infty}}}\|v(x,t)\| \leq M, \quad  \sup_{t\in \BR_+}\|v_x(0,t)\|<\infty.
\end{align}
Then, using \eqref{p22} (see, e.g., \cite[formula (4.14)]{FKRSsk}), we derive inequality \eqref{r9} for $R(t,z)=R(0,t,z)$ in a domain
\begin{align}& \label{r17}
\cld=\{z: \quad \Im(z) > M_1, \quad \Re(z) > M_2\},
\end{align}
and with some $\wh M>0$, $M_1>M$, $M_2>0$.  After that, the same considerations as in the scheme of the
proof of Theorem \ref{TmQP} yield our next theorem.
\begin{Tm}\label{TmQPsk} Let $v$  satisfy the conditions of Theorem    \ref{evol} $($for each $a\in \BR_+)$
and inequalities \eqref{r16}. Then, in some domain \eqref{r17}
we have the equality \\ $\vp(0,z) = -\lim_{t \to \infty} R_{22}(t,z)^{-1}R_{21}(t,z)$.
\end{Tm}
{\bf 3.} Finally, let us consider SGE \eqref{SG1}, where $v=-\psi_x$ is bounded in $\Om_{\infty}$:
\begin{align}& \label{r12}
\sup_{(x,t)\in \Om_{\infty}}| \psi_x(x,t) | \leq M.
\end{align}
The next proposition follows (see \cite[Corollary 6.25]{SaSaR}) from Theorem~\ref{SGevolSG}.
\begin{Pn}  \label{SGCyBVP} Let  $\psi$ given in $\Om_{\infty}$ satisfy  inequality \eqref{r12}.
Assume also that 
the conditions of  Theorem~\ref{SGevolSG} and  initial-boundary conditions \eqref{SG5}  hold
for each $\Om_a$ $(a\in \BR_+)$.
Then, for values
of $z\in \BC_M$, such that the inequalities
\begin{equation} \label{r13}
\big(\cos (2h_2(t))-\ve(z) \big)\Im(z) \geq |\Re(z)\sin ( 2h_2(t))|
\end{equation}
are valid for some $\ve(z)>0$ and for all $t \geq 0$, we have
\begin{equation} \label{r14}
\vp(0,z)=-\lim_{t \to \infty}R_{21}(t,z)/R_{22}(t,z).
\end{equation}
\end{Pn}
\subsection{Unbounded solutions of SGE}\label{unbound}
Theorem \ref{SGWsg}, Corollary \ref{SGMCy} and Proposition \ref{SGCyBVP} allow us to construct wide classes
of unbounded (in $\Om_{\infty}$) solutions of SGE, that is, solutions which do not satisfy \eqref{r12}. The simplest
case is the case $h_2(t)\equiv 0$ (see \cite[Proposition 6.27]{SaSaR}).
In this case, \eqref{r13} holds for all $z\in \BC_M$ ($M>0$) and, if other conditions of  Proposition \ref{SGCyBVP}  are valid,
\eqref{r14}
yields $\vp(0,z)\equiv 0$. Since $\M(0)\equiv 0$, we derive $h_1\equiv 0$. In other words, any solution of SGE with $h_2(t)\equiv 0$ and $h_1(x) \not\equiv 0$
does not satisfy conditions of  Proposition \ref{SGCyBVP}.  Solutions of initial-boundary value problems for SGE
may be constructed using Theorem \ref{SGWsg} and Corollary \ref{SGMCy}.
\begin{Pn}\label{Unb} Assume  that
$h_1(x)=\ov{h_1(x)}\not\equiv 0$ is three times  differentiable for $x \geq 0$, that 
\[h_1^{\prime}, \, h_1^{\prime \prime}, \, 
h_1^{\prime \prime \prime}
\in L^1(\BR_+), \quad h_1(0)=0,
\]
and that  $h_2\equiv 0$. Then one can use
the procedure given in Theorem \ref{SGWsg} in order to construct a solution
$\psi$ of the initial-boundary value problem  \eqref{SG1}, \eqref{SG5} for SGE,
and  the absolute value of the derivative  $\frac{\prt}{\prt x}\psi$ $($for the constructed $\psi)$ is always unbounded in the quarter-plane $\Om_{\infty}$.
\end{Pn}
Indeed, the solution constructed in Proposition \ref{Unb} satisfies  the conditions of  Proposition \ref{SGCyBVP}
excluding, possibly, \eqref{r12}, and so \eqref{r12} does not hold.

Some classes of unbounded solutions of the KdV equation with the minus sign before the dispersion term
are constructed in \cite{KdV} using low energy asymptotics of the Weyl functions.

\subsection{Reduction of the initial-boundary conditions in a semistrip}
{\bf 1.}  In this paragraph we present some results from \cite[Section 4]{ALSarx14} on the $N$-wave equation \eqref{NW1}
in the semistrip $\Om_a$. By $\vp(t,z)$ we denote evolution of the $GW$-function corresponding to this equation or, more precisely,    
$\vp(t,z)$ are
$GW$-functions of the systems $y_x=G(x,t,z)y$, where $G$ is expressed via $\vr$ using  \eqref{NW4} and \eqref{NW4'}.
The matrix function $R(t,z)=R(0,t,z)$ is determined via \eqref{p22}, where $F$ (similar to $G$) is given by \eqref{NW4} and \eqref{NW4'}.
First, we express  normalized (by \eqref{NW8}) $GW$-functions $\vp(t,z)$ via $\vp(0,z)$ and $R(t,z)$ (in other words, we derive the evolution
of the normalized $GW$-functions). Next, we consider the subcase 
when the entries of $\wh D$ are ordered in the same way as the entries of $D$:
\begin{align}  \label{NWsub3} &
\wh D=\di \{ \wh  d_1, \, \wh  d_2, \ldots,
\wh  d_m \}, \quad \wh  d_1> \wh d_2>\ldots> \wh d_m>0.
\end{align} 
\begin{Tm}
\label{TmNWevol}
Let $\vr=\vr^*$ satisfy the $N$-wave equation \eqref{NW1}, where $D$ and $\wh D$ have the form \eqref{NW2} and \eqref{NW3},
respectively. 
Assume that  $\vr(x,t)$ is uniformly bounded and continuously differentiable on $\Om_a$.

Then the matrix functions
\begin{align}\nn
 \psi_k(t,z):=&
 \begin{bmatrix} I_k & 0\end{bmatrix}
R(t,z)   \varphi (0,z)
   \left[
        \begin{array}{c}
              0       \\
              I_{m-k}
        \end{array}
   \right]
\\  \label{NW21} & \qquad \times   
   \biggl(
   \begin{bmatrix} 0 & I_{m-k}\end{bmatrix}
     R(t,z) \varphi (0,z)
      \left[
        \begin{array}{c}
              0       \\
              I_{m-k}
        \end{array}
      \right]
   \biggr)^{-1} ,
\end{align}
where $\vp(0,z)$ is the normalized $GW$-function of the system $y_x=G(x,0,z)$,  are well-defined for $1\leq k <m$. The  normalized $GW$-functions $\vp(t,z)$   are given 
$($in $\BC^{-}_M$ for some $M>0)$ by the formula
\begin{align}& \label{NW22}
      \left\{
            \varphi_{i, k+1}(t,z)
      \right\}_{i=1}^k
   =   \psi_k (t,z)\begin{bmatrix} 1\\ 0 \\ \ldots \\0
   \end{bmatrix}    \quad (z \in \BC^{-}_M ),
\end{align}
and by the normalization conditions \eqref{NW8}.
\end{Tm}
From Theorems \ref{NWTmUniq} and \ref{TmNWevol} follows the next result.
\begin{Tm}\label{TmNWnm} For the case where the entries of the matrix $\wh D$ in \eqref{NW1} are positive and ordered as in \eqref{NWsub3},
there is no more than one uniformly bounded and continuously 
differentiable on $\Om_a$ solution $\vr=\vr^*$  $($of the $N$-wave equation \eqref{NW1}$)$, having
the initial values $\vr(x,0)$  such
that
 $\vp(0,z)$ is bounded and \eqref{NW8!} holds. That is, there is no more than one
solution of the corresponding initial value problem.
\end{Tm}
Taking into account the results on fundamental solutions from \cite{BC0}, the requirements on $\vp(0,z)$ given in Theorem \ref{TmNWnm}   are reformulated below in terms of the sufficient
requirements on the initial condition
\begin{align}& \label{NW23}
\vr(x,0)=\rho(x)=\rho(x)^* \quad (0\leq x<\infty).
\end{align}
\begin{Pn} \label{NWPnNorm} Suppose that the initial condition $\rho(x)$ is absolutely continuous on $[0, \infty)$
and  the entries of $\rho(x)$ and $\rho^{\prime}(x)$ belong
$L^1(\BR_+)$. Then,  the normalized $GW$-function $\vp(0,z)$ of the system 
\begin{align}& \label{NW37}
y^{\prime}(x,z)=(\I z D-\zeta(x))y(x,z) \quad (x \geq 0), \quad \zeta=[D,  \rho]
\end{align}
is analytic and  bounded in $\BC^{-}_M$ $($for some $M>0)$  and \eqref{NW8!} is valid.
\end{Pn}
Theorem \ref{TmNWnm} is proved in \cite[Section 4]{ALSarx14} using the fact that 
\begin{align}& \label{NW22!}
\sup_{t\in [0,a), \, \Im(z)<-M}\|R(t,z)\vp(0,z)\exp\{-\I z t \wh D\}\|<\infty ,
\end{align}
which yields (see \cite[pp. 108, 109]{SaSaR}) that $\vp(0,z)$ determines not only the initial but also the boundary values of $\vr$.
In particular, when $\vp(0,z)$ satisfies the conditions of Theorem \ref{NWInvPr}, then, according to Remark \ref{RkPrM}, we have a procedure
to recover the initial and boundary conditions from $\vp(0,z)$. More precisely, we have 
\begin{align}& \label{exi1}
\zeta(x,0) =\M(D,\vp(0,z)),  \quad \wh \zeta(0,t) =\M(\wh D,\vp(0,z)),
\end{align}
and, using \eqref{NW4'}, we easily obtain $\vr(x,0)$ and $\vr(0,t)$. Indeed, without loss of generality we may assume that all the entries
of $\vr$ on the main diagonal are equal to zero and recover $\vr=\{\vr_{ik}\}_{i,k=1}^m$ from $\zeta$ or $\wh \zeta$ via formulas
\begin{align}& \nn
\vr_{ii} \equiv 0, \quad \vr_{ik}=\zeta_{ik}/(d_i - d_k)=\wh \zeta_{ik}/(\wh d_i - \wh  d_k) \,\, {\mathrm{for}} \,\, i\not= k.
\end{align}

{\bf 2.} In this and the next paragraph we deal with the cases of  dNLS \eqref{p3} with quasi-analytic boundary or initial, respectively, conditions,
 (i.e., the corresponding boundary or initial value functions belong to quasi-analytic classes). The definition of  quasi-analytic classes $C([0,\ell); \, \wt M)$ is given in Appendix \ref{AppB}. 

Further we present some results from \cite[Section 3]{ALSarx14}.
Let us consider  $m_1 \times m_2$ matrix functions $v(x,t)$,
which are continuously  differentiable and 
such that $v_{xx}$ exists
on the semi-strip $\Om_a$. 
Moreover, we require some smoothness of $v$ in the neighborhood of the point $(0,0)$. More precisely,
we require that for each $k$ there is a value $\ve_k=\ve_k(v)>0$ such that $v$
is $k$ times continuously differentiable with respect to $x$ in the square 
\begin{align}& \label{B3}
\cld({\ve_k})=\{(x,\, t): \quad  0\leq x\leq \ve_k, \quad 0\leq t\leq \ve_k\}, \quad \cld({\ve_k})\subset \Om_a.
\end{align}
The class of such functions $v(x,t)$ is denoted by $C_{\ve}(\Om_a)$.

\begin{Pn} \label{Cy3.4} Assume that  $v(x,t)\in C_{\ve}(\Om_a)$ satisfies  the dNLS equation \eqref{p3} on $\Om_a$
and that   a matrix function $v(0,t)$
or $v_x(0,t)$ is quasi-analytic. 

Then a matrix function $v(0,t)$ or $v_x(0,t)$, respectively, is uniquely determined
by the initial condition 
\begin{align}& \label{B4}
v(x,0)=h(x).
\end{align}
\end{Pn}
Proposition \ref{Cy3.4} is proved by presenting formulas to recover the derivatives  $\left(\frac{\p^k}{\p t^k}v\right)(0,0)\,$ and $\, \left(\frac{\p^k}{\p t^k}v_x\right)(0,0)\,$
$\,(k \geq 0)$ from $\, h(x) \,$ (see the proof of   \cite[Proposition 3.2]{ALSarx14}).

Using Proposition \ref{PnWsead}, Theorems \ref{PrIP} and \ref{TmevNLS} and Proposition \ref{Cy3.4} we derive the next theorem.
\begin{Tm} Assume that  $v\in C_{\ve}(\Om_a)$ satisfies  the dNLS equation \eqref{p3} on $\Om_a$, that \eqref{p29} holds
 and that the functions $v(0,t)$ and $v_x(0,t)$ $($boundary values$)$ belong to some quasi-analytic classes $ C([0,a); \, \wt M)$ and $C([0,a); \, \wt M^+)$,
 respectively. Then $v$ is uniquely  determined by the initial condition \eqref{B4}.
\end{Tm}
{\bf 3.}  Using Proposition \ref{PnWsead} and Theorems \ref{PrIP} and \ref{TmevNLS} one can also derive Theorem 4.3 from \cite{Sa14b},
which we formulate in this paragraph. The proof is based on the formulas to recover the derivatives  $\left(\frac{\p^k}{\p x^k}v\right)(0,0)$ 
$(k \geq 0)$ from the boundary conditions
\begin{align}& \label{B5}
v(0,t)=h_0(t), \quad v_x(0,t)=h_1(t) \quad (0\leq t<a).
\end{align}
\begin{Tm} Assume that  $v\in C_{\ve}(\Om_a)$ satisfies  the dNLS equation \eqref{p3} on $\Om_a$, that \eqref{p29} holds
 and that the  initial value function  $v(x,0)$  belongs to some quasi-analytic class $ C([0,\infty); \, \wt M)$.
 Then $v$ is uniquely determined by the boundary conditions \eqref{B5}.
\end{Tm}

\section{Dynamical Dirac system and response function} \label{Dyn}
\setcounter{equation}{0} 
\subsection{Introduction}
In this section we give several statements  from our paper \cite{SaA15Dyn}, which appeared in arXiv in 2015.
Classical Dirac systems \eqref{p1} are also called spectral Dirac systems and various new results
on inverse problems for these systems appeared quite recently (see some references in Subsection \ref{Mais}). 
 At the same time, a great and growing interest in
dynamical systems and control methods is reflected in the active study of the dynamical inverse problems
and, in particular, in the study of the inverse problems for dynamical Schr\"odinger and Dirac systems \cite{AMRa, AMRy, Bel2, BeMi, GladMo, MeOR}
(see also the references therein).

Dynamical Dirac system
(Dirac system in 
the {\it time-domain setup}) was studied in the important recent paper \cite{BeMi}.  The dynamical Dirac system considered in \cite{BeMi} is an evolution system of hyperbolic type
and has the following  form:
\begin{align}\label{td1.3}
& \I Y_t+JY_x+\clv Y=0 \quad (x>0, \quad t>0); \\ & \label{td1.3'}
Y=\begin{bmatrix}Y_1 \\ Y_2 \end{bmatrix}, \quad J=\begin{bmatrix} 0 & 1 \\  -1 & 0 \end{bmatrix}, \quad \clv=\begin{bmatrix} p & q \\  q & -p \end{bmatrix}, \quad Y_t:=\frac{\p Y}{\p t},
\end{align}
where $p=p(x)$ and $q=q(x)$ are real-valued functions of $x$, and initial-boundary conditions are given by the equalities
\begin{align}\label{td1.4}
Y(x,0)=0, \,\, x \geq 0; \quad Y_1(0,t)=f(t), \,\, t \geq 0.
\end{align}
Here $f$ is a complex-valued function (called boundary control by \cite{BeMi}) and the input-output map (response operator) $R: \, Y_1(0,\cdot) \rightarrow Y_2(0,\cdot)$
is of the convolution form $Rf=\I f +r*f$. Thus,  $R$ and $r$ denote in this section the response operator and response function, respectively. ($\wh Y$ stands in this section
for the Fourier transformation of $Y$, see \eqref{td2.13}.)

The inverse problem consists in recovery of the potential $\clv$ from the {\it response function} $r.$
This  inverse problem was considered in  \cite{BeMi} using boundary control methods.

We note that recent results on dynamical Schr\"odinger and Dirac equations are based on several earlier works.
In his paper \cite{Blag} from 1971, A.S.~Bla-gove\v{s}\v{c}enskii considered dynamical system
\begin{align} \label{h1}&
Y_{tt}-Y_{xx}+ \clq (x)Y_x=0
\end{align}
with boundary control $Y(0,t)=f(t)$, and solved inverse problem to recover $\clq$ from $f$. A.S. Blagove\v{s}\v{c}enskii
established important connections between his problem and spectral theory of string equations.
This work was developed further in \cite{Bel1} (see also references therein), where response operator appears
in  inverse problem. Finally, the inverse problem to recover the matrix potential $\clq(x)$ of the dynamical Schr\"odinger equation
\begin{align} \label{h2}&
Y_{tt}-Y_{xx}+ \clq (x)Y=0
\end{align} 
from the response function was considered in \cite{ABI} by  S. Avdonin, M. Belishev, and S. Ivanov. 

Similar to the case of the spectral Dirac and Schr\"odinger equations, the dynamical Dirac  equation is
a more general object than the dynamical Schr\"odinger equation. More precisely, setting
in \eqref{td1.3} 
\begin{align} \label{h3}&
p(x)=0, \quad q(x)=g_x(x)/g(x), \quad {\mathrm{where}} \quad g_{xx}(x)=  \clq (x)g(x),
\end{align} 
and rewriting \eqref{td1.3}, \eqref{td1.3'} in the form
\begin{align} \label{h4}&
(Y_1)_t=\I\big((Y_2)_x+qY_2\big), \quad (Y_2)_t=\I\big(-(Y_1)_x+qY_1\big),
\end{align} 
we obtain a dynamical Schr\"odinger equation
\begin{align} \nonumber
(Y_1)_{tt}&=(Y_1)_{xx}-(q_x+q^2)Y_1=(Y_1)_{xx}-\big(g_{xx}/g\big)Y_1
\\  \label{h5}&
=(Y_1)_{xx}-\clq Y_1.
\end{align} 
For interesting applications of the interconnections between  spectral Dirac and Schr\"odinger equations
see, for instance, the papers \cite{BolGe, CGR, GeGo, GeSS, EGNT, EGNST} and references therein.
\subsection{Preliminaries and estimates} \label{Est}
According to \cite[Theorem 1]{BeMi}, in the case where
 $p, \, q, \, f$ (in \eqref{td1.3'} and in \eqref{td1.4}) are continuously differentiable and $f(0)=f^{\prime}(0)=0$, there
is a unique classical solution $Y$ of  \eqref{td1.3}, \eqref{td1.4}  and this solution
admits representation
\begin{align}
\label{td2.8}&
 Y=Y_*^f+w^f; \qquad Y_*^f(x,t)=f(t-x)\begin{bmatrix} 1 \\ \I \end{bmatrix},  \quad f(t)=0 \quad (t < 0);\\
\label{td2.8'}& w^f(x,t)= \int_x^t f(t-s)\kappa(x,s)ds  \quad (t\geq x\geq 0), 
\\ \label{td2.8''}&
 w^f(x,t)=0  \quad (x> t\geq 0),
\end{align}
where $\kappa(x,s)$ $(x \leq s)$ is continuously differentiable. In particular, formulas \eqref{td2.8} and \eqref{td2.8''} yield:
\begin{align} \label{td2.f}&
Y(x,t)=0 \,\, {\mathrm{for}} \,\, 0 \leq t<x \quad({\mathrm{finiteness \, of \, the \, domain \, of \, influence }}).
\end{align}
Representation \eqref{td2.8}--\eqref{td2.8''} is proved in \cite{BeMi} using Duhamel formula. 

Let us also assume that $\clv$, $f$ and $f^{\prime}$ are bounded:
\begin{align} \label{td2.11}&
\sup_{x >0}\|\clv(x) \|<M_1, \quad \sup_{t>0}\left\|f(t)\begin{bmatrix} 1 \\ \I \end{bmatrix}\right\|<c_0 ,
\quad \sup_{t>0}\left\| f^{\prime}(t)\begin{bmatrix} 1 \\ \I \end{bmatrix}\right\|< \wt c_0.
\end{align}
The following estimates are proved in \cite[Section 2]{SaA15Dyn}.
\begin{Pn} \label{PnEst} Let $p, \, q, \, f$ be continuously differentiable  and let equalities  $f(t)=f^{\prime}(t)=0$ hold for
$t \leq 0$. Assume that \eqref{td2.11} is valid. Then the solution $Y$ of the dynamical Dirac system \eqref{td1.3},
such that \eqref{td1.4} and \eqref{td2.f} are valid, satisfies the following inequalities
\begin{align} \label{td2.!}&
\| Y(x,t)\| \leq c_0\E^{Mt}, \quad \| Y_t(x,t)\| \leq \wt c_0\E^{Mt}, \quad \| Y_x(x,t)\| \leq M_2 \E^{Mt},
\end{align}
where $x\geq 0$, $t\geq 0$,  $M_2>0$ is some constant, and $M=2 \sqrt{2} M_1$.

The kernel $\kappa$ of the integral operator in \eqref{td2.8'} satifies the inequality
\begin{align} \label{tdp3}&
\| \kappa(x,t)\| \leq M\E^{Mt}.
\end{align}
\end{Pn} 
In view of \eqref{td2.!}, we can apply to $Y$ the transformation: 
\begin{align} \label{td2.13}&
\wh Y(x,z)=\int_{0}^{\infty}\E^{\I z t} Y(x,t)dt, \quad z \in \BC_M,
\end{align}
and $\wh Y$ stands in this section for the Fourier transformation of $Y$ $($which is taken, for the sake of convenience, for the fixed values  $x,z)$. 
Moreover, the same transformation can be applied to $Y_t$ and $Y_x$, and we have
\begin{align} \label{td2.13'}&
\I \int_{0}^{\infty}\E^{\I z t} Y_t(x,t)dt =z\wh Y(x,z),
\\ \label{td2.13''}&
\int_{0}^{\infty}\E^{\I z t} Y_x(x,t)dt =\frac{d}{d x}\wh Y(x,z)=:\wh Y^{\prime}(x,z),\quad z \in \BC_M.
\end{align} 
Now, applying the Fourier transformation to the dynamical Dirac system \eqref{td1.3}, we derive
\begin{align} \label{td2.16}&
z\wh Y(x,z) +J \wh Y^{\prime}(x,z) +\clv(x)\wh Y(x,z)=0.
\end{align}
Note that, according to \cite{BeMi}, the {\it response function} $r$ is given by
\begin{align}& \label{tdp4-}
r(t)=\kappa_2(0,t).
\end{align}
Indeed, for $r$ of the form \eqref{tdp4-}, using  \eqref{td1.4}, \eqref{td2.8} and \eqref{td2.8'} we obtain
$Y_1(0,t)=f(t)$ and
\begin{equation}\label{tdp4}
 Y_2(0,t)=\I f(t)+ \int_0^t r(t-s)f(s)ds=\I Y_1(0,t)+\int_0^t r(t-s)Y_1(0,s)ds.
\end{equation}
\subsection{Response and Weyl functions} \label{Resp}
In this subsection we always assume that the conditions of Proposition \ref{PnEst} are fulfilled.
However,  it  seems possible (and would be  interesting) to modify representation \eqref{td2.8}-\eqref{td2.8''}
for the case of matrix functions $p$ and $q$ and for weaker smoothness conditions, in which case the solution of the
inverse problem in a much more general situation will follow.

The equivalence transformation between spectral Dirac systems in the form \eqref{td2.16}, where $J$ and $\clv$ are given in \eqref{td1.3'}, and in the form  \eqref{p1}, \eqref{p2} ($m_1=m_2=1$)
is given by the relations
\begin{align} \label{td2.19}&
y=K\wh Y, \quad v=\I q - p, \quad K= \frac{1}{\sqrt{2}} \begin{bmatrix} \I & 1 \\ - \I & 1 \end{bmatrix}.
\end{align}
Since, according to \cite[Section 3]{SaA15Dyn}, $\wh u(x,z)$ ($z\in \BC_M$) is a Weyl solution of  \eqref{td2.16} (i.e., $\wh u \in L^2_2(0,\infty)$),
we see that  $y=K\wh Y$ is a Weyl solution of \eqref{p1}, \eqref{p2}, where $m_1=m_2=1$. Hence, in view of \eqref{tdp4}, definition \eqref{p7} yields \cite[Proposition 3.1]{SaA15Dyn}:
\begin{Pn}\label{tdPn1} The response function $r(t)$ of the dynamical Dirac system \eqref{td1.3} is connected
with the Weyl function $\vp(z)$ of the corresponding spectral Dirac system \eqref{p1} $($where
$m_1=m_2=1)$ via the equality 
\begin{align}& \label{td2.22}
\vp(z)=\wh r(z)/(\wh r(z)+2\I), \quad z\in \BC_M.
\end{align}
\end{Pn}
\begin{Rk}\label{tdRkIP} Proposition \ref{tdPn1} jointly with Theorem \ref{PrIP} and the equality $v=\I q - p$ in the equivalence transformation \eqref{td2.19}
$($which yields $p=-\Re(v)$,  $q~=~\Im(v))$ provide a procedure to solve the inverse problem of the recovery of the potential $\clv$ of a dynamical Dirac system
from its response function.
\end{Rk}
Since $m_1=m_2$, it is possible to introduce a slightly different class of Weyl functions  (that is, Weyl functions  $\vp_H$) via the inequality:
\begin{align} \label{td48}&
\int_0^\infty \left[ \begin{array}{lr} I_k &  \I \varphi_H (z )^*
\end{array} \right]
  \Theta  u(x, z )^*
 u(x, z )\Theta ^*
 \left[ \begin{array}{c}
I_k \\ - \I \varphi_H (z ) \end{array} \right] dx < \infty , 
\\ \label{th1}&
{\Theta}:=   \frac{1}{\sqrt{2}}       \left[
\begin{array}{cc} I_k &
-I_{k} \\ I_{k} & I_k
\end{array}
\right],  \quad z \in \BC_+,
\end{align}
where $k=m_1=m_2$. Comparing \eqref{td48} with \eqref{p7} we see that $\vp_H$ is  a linear fractional transformation of $\vp$.
It is convenient that Weyl functions $\vp_H$ belong to Herglotz class (i.e., $\I(\vp_H(z)^*-\vp_H(z))\geq 0$).
In our case we have $k=1$ and the connection between $\vp_H$ and $\wh r$ is simpler than \eqref{td2.22}. Namely, we have
\begin{align}& \label{td2.23}
\vp_H(z)=\wh r(z)+\I .
\end{align}
\begin{Rk} \label{tdRkHer} Clearly, we can recover $v$ from $\vp_H$ after an easy transformation which maps $\vp_H$ into $\vp$.
However, there is also an independent procedure in \cite{SaA4} $($see \cite{SaA15Dyn, SaSaR} as well$)$ to recover $v$ from $\vp_H$.
Instead of the structured operators $S_{ l}$ given by \eqref{p15} and \eqref{p13}, convolution operators $\wt S_{ l}:$
\begin{align}& \label{tds1}
\wt S_{ l}=\frac{d}{dx}\int_0^{ l} \wt s(x-t) \, \cdot \, dt, \qquad \wt s(x)=- \wt s(-x)^* \quad (x >0),
\\ & \label{tds2}
\wt s(x)^*:=  \frac{d}{dx}\left( \frac{\I}{4 \pi}
\E^{ \eta x}
{\mathrm{l.i.m.}}_{a\to \infty} \int_{- a}^{a}\E^{-\I \xi x}
(\xi +\I \eta)^{-2} \varphi_H (\xi +\I \eta)d \xi \right)
\end{align}
are used for this purpose.  Formulas \eqref{td2.23} and \eqref{tds2} imply that
$r(t)=2\I \ov{\wt s^{\prime}(t)}$. Moreover, it is shown in \cite{SaA15Dyn} that  the response function $r$ coincides
with a so called $A$-amplitude from \cite{GeSi} $($more precisely, with the analogue  of $A$-amplitude for the corresponding spectral Dirac system$)$. We note that in M.G. Krein's terminology
$\wt s^{\prime}$ is called the accelerant.
\end{Rk}
Interconnections between response functions, Weyl functions and $A$-amplitudes for Schr\"odinger equations
are discussed in the interesting paper \cite{AMRy}.

Some explicit solutions of the inverse problem for dynamical Dirac system are obtained in \cite[Section 5]{SaA15Dyn}:
\begin{Tm}\label{TmEDIP} \cite{SaA15Dyn} Let $r(t)$ be the response function of a dynamical Dirac system  
and assume that $r(t)$ admits representation  
\begin{align}\nn &
r(t)=-2\I{\vt_2^*}\E^{-\I t \a}\vt_1,
\end{align}
where the $n \times n$ $(n \in \BN)$ matrix
$\a$ and the column vectors   $\vt_i \in \BC^N$ $(i=1,2)$ satisfy the identity
$$\a-\a^*=-\I\big(\vt_1+\vt_2\big)\big(\vt_1+\vt_2\big)^*.$$
Then the potential $\clv$ of this dynamical Dirac system
is given $($in terms of $\a$ and $\vt_i)$ by the equalities
\begin{align}\nn
& \clv=\begin{bmatrix} p & q \\  q & -p \end{bmatrix},  \quad 
p=-\Re(v),  \quad q=\Im(v);
\\ \nn &
v(x)=-2\I  \vt_1^*\E^{\I xA^*}S(x)^{-1}\E^{\I xA}\vt_2, \quad A:=\a+\I \vt_1
\big(\vt_1+\vt_2\big)^*;
\\ \nn &
S(x)=I_n+\int_0^x \Lam(t)\Lam(t)^*dt, \quad
\Lam(t)
=\begin{bmatrix}\E^{-\I tA}\vt_1 \quad \E^{\I tA}\vt_2\end{bmatrix}.
\end{align}

\end{Tm}

\bigskip 
\noindent{\bf Acknowledgments.}
 {This research  was supported by the
Austrian Science Fund (FWF) under Grant  No. P24301.}

\appendix

\section{Quasi-analytic functions and matrix functions}\label{AppB}
\setcounter{equation}{0} 
The class $C\big(\{\wt M_k\}\big)$  
consists of all infinitely differentiable on $[0,  \ell)$ scalar functions $f$ such that for some 
$c(f) \geq 0$ and  for fixed constants $\wt M_k >0 $ $(k \geq 0)$ we have
\begin{align}& \label{B1}
\left| \frac{d^k f}{dt^k}(t)\right|\leq (c(f))^{k+1}\wt M_k \quad {\mathrm{for \, all}} \quad t\in [0,  \ell), \quad 0< \ell\leq \infty.
\end{align}
Recall that $C\big(\{\wt M_k\}\big)$    is called quasi-analytic if 
for the functions $f$ from this class and for any $0\leq t < \ell$ the equalities $\frac{d^k f}{dt^k}(t)=0$ $(k\geq 0)$ yield $f\equiv 0$.
According to the famous Denjoy--Carleman theorem, the equality
\begin{align}& \label{B2}
\sum_{n=1}^{\infty}\frac{1}{L_n}=\infty, \quad L_n:=\inf_{k\geq n}\wt M_k^{1/k}
\end{align}
implies that  the class $C\big(\{\wt M_k\}\big)$ is quasi-analytic. 

In the case of matrix functions $\phi(t)$ (e.g, in the case of $\phi(t) = v(0,t)$) we say that $\phi$ is quasi-analytic if the
entries of $\phi$ are quasi-analytic and we say that $\phi \in C([0, \ell); \, \wt M)$, where $\wt M=\{\wt M_k(i,j)\} $,
if the entries $\phi_{ij}$ of $\phi$ belong to  quasi-analytic classes $C(\{\wt M_k(i,j)\} )$.

The recovery of a  function belonging to some quasi-analytic class $C\big(\{\wt M_k\}\big)$ from its Taylor coefficients (and from the sequence $\{\wt M_k\}$)
is discussed
in important works   \cite{Bang, Khr} and \cite[Section III.8]{Beur}  but many interesting problems are still open.




\begin{flushright}
Alexander Sakhnovich,\\
Vienna
University
of
Technology,
Austria, \\
e-mail: {\tt oleksandr.sakhnovych@tuwien.ac.at}

\end{flushright}

\end{document}